\documentclass[12pt]{amsart}
\usepackage[T1]{fontenc}
\usepackage{amssymb}
\newtheorem{theorem}{Theorem}

\author{Sung Soo Kim}
\address{Sung Soo Kim \\
 Hanyang University, Ansan, Kyunggi 425-791, Korea}
\email{kim5466@hanyang.ac.kr}
\author{Szymon Plewik}
\address{Szymon Plewik\\Institute of Mathematics,
University of Silesia, ul. Ban\-ko\-wa 14, 40-007 Katowice}
\email{plewik@ux2.math.us.edu.pl}

\begin{document}

 \title{DISCONTINUITY  AND INVOLUTIONS ON COUNTABLE SETS}

\begin{abstract}  
For any infinite subset $X$ of the rationals and a subset
$F \subseteq X$ which has no isolated points in $X$ we  construct
 a function $f: X \to X$ such that  $f(f(x))=x$ for each $x\in X$ and  $F
$ is the set of discontinuity points of  $f$.
\end{abstract} 

\maketitle
 In the literature one  finds a few
algorithms that can produce any given subset of the rationals as
the set of discontinuity points of a function.  Probably Wac\l aw
Sierpi\'{n}ski [2] was the first to publish the algorithm, of the
kind that is best known. In [1] this algorithm was reduced to
following: let $X = A \cup B $ be a topological space, where sets
$A$ and $B$ are dense and disjoint; assume that $Y = \{ 0\} \cup
\{ {1\over{n}}: n=1,2,\ldots \} \cup \{ {-1\over{n}}:
n=1,2,\ldots \}$; suppose that  $X \setminus C$ is the
intersection of a decreasing sequence of open sets $F_n
\subseteq X$ with $F_1 = X$; if $x\in X \setminus C,$ then put
$f(x) =0$; if $x\in A\cap  F_n \setminus F_{n+1},$
 then put $f(x) ={1\over{n}}$; if $x\in B\cap  F_n \setminus
F_{n+1},$ then put $f(x) ={-1\over{n}}$; the set $C$ is the set of
discontinuity points of the defined function $f:  X \to Y$.  In
this note we are suggesting an algorithm that works with involutions.

 Let us assume that $F$ and $Q$ are disjoint
subsets of the rationals.

\begin{theorem} If $F$ is infinite and $F$ has no
isolated point in $Q\cup F$, then there is a bijection $f: Q\cup
F \to Q \cup F$ such that:  $Q$ is the set of continuity points
of $f$; $f$ is the identity on $Q$;  for any $x \in F$ we have
$f(x)\not= x$ and $f(f(x)) =x$.
\end{theorem}

 \begin{proof} Enumerate all points of $Q$ as  a sequence
$y_0, y_1,
  \ldots $; enumerate all points of $F$ as a sequence $x_0, x_1,
  \ldots $; choose an irrational number $g$ such that   $F\cap
(-\infty,
  g)$  is empty or infinite, and   $F\cap (g,
  +\infty)$  is empty or infinite; put $G_0 = \{ (-\infty, g),
  (g,+\infty)\}.$

  Take $x_0$ and choose $f(x_0) \in F \cap A$ such that
$f(x_0)\not= x_{0}
  \in A \in G_0$.   Put $f(f(x_0)) =x_0$ and $F_0 = \{ -\infty ,
  +\infty , g, x_0, f(x_0), y_0 \} .$ Let $G_1$ be a family of all
open
  intervals with endpoints which are succeeding points of $F_0$.
  Suppose that the set $F_n$ has been defined and let $G_{n+1}$
  be consisted of all intervals with endpoints which are succeeding
points of
  $F_n$. Let  $x_{k_n} \in F \setminus F_n$  be the  point with the
least
  possible  index such that  $f(x_{k_n})$ has  not been defined,
  but $f(x_i)$ has
  been defined for any $i<k_n$. Choose $f(x_{k_n}) \in F \cap
A\setminus F_n$
  such that $f(x_{k_n})\not= x_{k_n} \in A \in G_j$, where $j \leq
n+1$ is the
  greatest natural number for which a suitable $f(x_{k_n})$ could be
chosen.
  Put $ f( f ( x_{k_n} ) ) = x_{k_n} $ and $ F_{n+1} = F_n \cup \{
x_{k_n}, f(x_{k_n}),
  y_{n+1} \} .$ The bijection $f$ also requires that we set  $f(y_n)
= y_n$ for every  $n$. The
combinatorial properties of $f$ follow directly from the
definition. However, it remains to examine the continuity and
discontinuity of $f$.

Suppose $x \in F_m \cap F$ and $\{ a_0, a_1 , \ldots \} \subseteq
Q\cup F$ is a monotone sequence which converges to $x$.
Choose a
natural number $i\geq m$ such that for any $k\geq i$ there is
some   $I\in G_{k+1}$ and we have:  $x$ is an endpoint of $I$;
$a_n \in I$ for all but finite many $n$; $f(x )$ is not an
endpoint of $ I$. By the definition $f(a_n) \in I$ for all but
finite many $n$. It follows that $\lim_{n\to \infty } f(a_n) \not=
f(x )$. Therefore   $f$ is discontinuous at any point $x\in F$.

Note that if $y \in Q$ is an isolated point in $Q\cup F$, then
there is nothing to prove about the continuity of $f$ at $y$.
Suppose $y_m\in Q$ and $\{ a_0, a_1 , \ldots \} \subseteq Q\cup
F$ is a monotone sequence which converges to $y_m$. Then for
any
$k\geq m$ there is some   $I\in G_{k+1}$ and we have:  $y_m$ is
an
endpoint of $I$; $a_n \in I$  for all but finite many $n$. By the
definition $f(a_n) \in I$ for all but finite many $n$. It follows
that $\lim_{n\to \infty} a_n = y_m = f(y_m) = \lim_{n\to \infty }
f(a_n) $. Therefore  $f$ is continuous at any point $x\in Q$.
\end{proof}

{\bf References}

[1] S. S. Kim,  A Characterization of the Set of Points of Continuity
of a Real Function, {\it Amer. Math. Monthly}, 106 (1999), 258 -
259.

[2] W. Sierpi\'{n}ski, {it FUNKCJE  PRZEDSTAWIALNE
ANALITYCZNIE},
  Lw\'{o}w - Warszawa - Krak\'{o}w: Wydawnictwo
 Zak\l adu Narodowego Imienia Ossoli\'{n}skich (1925).

For related topic see:

P. R. Halmos, 
\textit{Permutations of sequences and the Schr¨oder-Bernstein theorem},
Proc. Amer. Math. Soc. 19 (1968), 509 - 510. MR0226590 (37 $\#$2179).

E. Hlawka, 
\textit{Folgen auf kompakten Räumen. II.} (German) 
Math. Nachr. 18 (1958) 188 - 202. MR0099556 (20 $\#$5995).
 
H. Niederreiter, 
\textit{A general rearrangement theorem for sequences}. 
Arch. Math. (Basel) 43 (1984), no. 6, 530--534. MR0775741 (86e:11061).

J. von Neumann, (1925)???

J.A. Yorke, 
\textit{Permutations and two sequences with the same cluster set}, 
Proc. Amer. Math. Soc. 20 (1969), 606. MR0235516 (38 $\#$3825).

\end{document}